\documentclass[12pt]{article}

\usepackage{amsfonts,amsmath,amsxtra}

=14.6mm  \relpenalty =1000 \binoppenalty=1000
\begin{document}
\centerline{\bf I. Mullayeva} \vspace{10mm} \centerline{\bf On the
weight structure of cyclic codes over $GF(q)$, $q>2$.}
\vspace{10mm}
\begin{abstract}

        The interrelation between the cyclic structure of an
    ideal, i.e., a cyclic code over Galois field $GF(q)$, $q>2$,
    and its classes of proportional elements is considered.
    This relation is used in order to define the code's weight
    structure. The equidistance conditions of irreducible
    nonprimitive codes over GF(q) are given. Besides that, the
    minimum distance for some class of nonprimitive cyclic
    codes is found.

\end{abstract}

The relation of proportionality for elements of algebra~$A_n$,
consisting of polynomials in $x$ over Galois field $GF(q)$, modulo
polynomial ${x^{n} - 1}$, is the equivalence relation ~\cite{VAN}.
Therefore $A_n$ falls into several disjoint subsets and every such
subset contains all elements which are proportional to each other.
These subsets will be called the classes of proportionality. Let
${z(x) \ne 0}$ be some vector of~$A_n$. If $\alpha_1=1, \alpha_2,
\ldots , \alpha_{q-1}$ are all different elements of the
multiplicative group~$GF(q)^*$ of the field $GF(q)$, then the
following $q - 1$ vectors
$$
\alpha_{1}z(x), \alpha_{2}z(x), \ldots , \alpha_{q-1}z(x)
\eqno(1)
$$
are some different and proportional to each other elements of
$A_n$. The set of vectors~$(1)$ is closed under the multiplication
by the elements of the group~$GF(q)^*$.  Hence the set $(1)$
represents some class of proportional elements, which will be
denoted by $P_{z(x)}$. Because of an arbitrary choice for $z(x)$,
every nonzero class  consists of $q-1$ elements of the form $(1)$.
Consequently $A_n$ contains $(q^{n}-1)/(q-1)$ different nonzero
classes. Evidently all elements of one class have the same period
~\cite{MAKS,NIL} or the same order \cite{LID}. Clearly, the
supporting sets \cite{MAKS} of vectors, entering into the same
class of proportionality, are similar too. Hence, the Hamming
weight is also the same for all vectors of one class. Thus, we can
say that any proportionality class ~$P_{z(x)}$, $P_{z(x)}\subset
A_n$, has its order, its supporting set and its Hamming weight.
Obviously, any proportionality class of $A_n$ is characterizied by
its unique monic polynomial.

Now consider an ideal $J$, $J\subset A_n$, i.e., some cyclic code
over ~$GF(q)$, having the following generator~$g(x)=(x^{n} -
1)/h(x)$~\cite{PET}, where $h(x)$~is some parity-check polynomial
of degree $m$, having some order $n$, $n=ord(h(x))$~\cite{LID}.
Below we suppose that $q>2$ and $gcd(n,q)=1$.

It's also known~\cite{NIL,MAKV,ELS} that any ideal is partitioned
into several disjoint subsets, that is cycles, under the
multiplication of ideal's vectors by $x$. On the other hand, some
ideal~$J$, as a subspace of $A_n$, consists of $(q^{m}-1)/(q-1)$
nonzero proportionality classes. Obviously, the existence of two
different partitions into some disjoint subsets of any ideal
assumes a certain dependence between proportionality classes and
cycles of ideal.

Further, any ideal  $J\subset A_n$ is the direct sum of minimal
ideals~ {\cite{MAKS,LID}}
$$
                        J= \sum_{i=1}^{t} J_i,                 \eqno
(2)
$$
where $J_i$ is some minimal ideal, having an irreducible
parity-check polynomial $h_{i}(x)$ of degree $m_i$ and of
order~$n_i$, $n_i=ord(h_{i}(x))$, ${1\le i\le t}$.  This implies
that the following polynomial

$$
               h(x) = \prod_{i=1}^t h_{i}(x),                  \eqno
(3)
$$
is the parity-check polynomial of $J$ and $n$, $n=lcm(n_1, n_2,
\ldots , n_t)$, is the order of $h(x)$~\cite{LID}. It should be
stressed that under the condition $gcd(n,q)=1$ the polynomial
$(3)$ has no repeated factors~\cite{LID}.

{\bf Remark $1$}. It is worth mentioning that the number $n$ can
be some number of either the primitive form~$n=q^{m}-1$  or of the
nonprimitive form $n\ne q^{m}-1$. In the first case, we have some
cyclic primitive code and the second case corresponds to a certain
cyclic nonprimitive code~\cite{BLH}. Let us stress that $n \ne
q^{m}-1$ if and only if $n_i \ne q^{m_i}-1$, $1\le i \le t$. But
if there is at least one primitive polynomial among the
polynomials $h_{i}(x)$, $1\le i\le t$, then $n=q^{m}-1$.

Furthermore, applying the theory of linear recurring
sequences~\cite{LID,ZIR} to elements of an ideal, we obtain that
every element of some ideal $J\subset A_n$ is characterized by its
unique minimal polynomial. Denote by ~$C$ the set of all elements
of~ $J$, having the same minimal polynomial $c(x)$. The set $C$ is
either some minimal ideal~$J_i$, $1\le i\le t$, or a certain
subset of all such elements of~$J$, whose characteristic
polynomial of the smallest degree coincides with~$c(x)$. In the
general case, the polynomial~$c(x)$ is equal to the product of
some $k$, $1\le k \le t$, polynomials from $t$ different prime
divisors of~$h(x)$. Thus,

$$
          c(x)= \prod_{j=1}^{k} h_{i_{j}}(x), \quad 1\le k\le t.
\eqno (4)
$$
This means that any element of $C$ has the same period or the same
order $n_c = ord(c(x))$, $n_c \le n$, $n_c | n$.

{\bf Lemma $1$.} {\it For the set} $C$, $C\subset J$, {\it having
some minimal polynomial}~$c(x)$ {\it in terms of}~$(4)$, {\it the
following equality takes place}
$$
           n_c \cdot s_c = R_{c}\cdot (q - 1),                    \eqno
(5)
$$
{\it where} $s_c$ {\it is the number of all cycles, and} $R_c$
{\it is the number of all proportionality classes of}~$C$.

{\bf Proof.} The set $C$ is closed under two different operations.
The first operation is the cyclic shift of vector and the second
one is the multiplication of vectors by elements of group
$GF(q)^*$. Hence equality~$(5)$ can be obtained by the counting of
the number of all elements, belonging to $C$, via the two
different ways. The lemma is proved.

{\bf Theorem $1$.} {\it Any cycle} $\{ z(x)\}$ {\it of ideal}
$(2)$, {\it having some period} $n_z$, $n_z | n$, {\it consists
of}~ $d_z$ {\it subsets}.  {\it The first element of each such
subset is proportional to} $z(x)$.  {\it Every such subset
contains} $r_z$ {\it nonproportional to each other vectors, i.e.,}
$n_z=r_z \cdot d_z$, $d_z|(q-1)$. {\it And the number} $r_z$ {\it
is the index of the subgroup, belonging to} $GF(q)^{*}$, {\it of
order~} $d_z$ {\it in the group of the roots of unity, having the
least possible order.}

{\bf Proof.} Let $r_z$, $1 \le r_z \le n_z $, is the smallest
natural number such that the following equality holds
$$
  x^{r_z}\cdot z(x)=\alpha \cdot z(x) mod(x^{n_z}-1),      \eqno
(6)
$$
where $\alpha$ is some element of $GF(q)^*$. Then the following
$r_z$ vectors of cycle $\{ z(x)\}$
$$
  z(x), x\cdot z(x),   ...   , x^{r^{z} - 1} \cdot z(x)
\eqno(7)
$$
are some non--proportional to each others vectors because,
assuming the inverse, we should be able to decrease~the number
$r_z$, but it is impossible.  Hence the set of elements $(7)$
belongs to the following $r_z$ classes of proportionality

$$
P_{z(x)}, P_{x z(x)}, \ldots    ,...,   P_{x^{r_{z}-1}z(x)}.
\eqno (8)
$$
Since $x^{n_z}z(x)=z(x)$ in the ring $A_n$ and also, considering~
$(6)$, we see that $P_{z(x)}=P_{x^{r_z}z(x)}=P_{x^{n_z}z(x)}$.
This means that the cycle~$\{ z(x)\}$ belongs to the classes~$(8)$
and every such class contains~$d_z$ vectors, $d_z=n_z/r_z$,
$1<d_z\le q-1$.  In terms of equality~ $(6)$ the following
different vectors $z(x), x^{r_z} z(x), x^{2r_z} z(x), \ldots ,
x^{(d_{z}-1)r_z}z(x)$ of class $P_{z(x)}$ can be represented as
${\alpha^{0}z(x), \alpha z(x),\ldots, \alpha ^{d_{z}-1}z(x)}$.
This implies that $\alpha^{0}=1, \alpha, \alpha^{2}, \ldots,
\alpha^{d_{z}-1}$ are the different elements of group $GF(q)^*$.
Since $x^{n_z}z(x)=x^{r_{z}d_{z}}z(x) =\alpha^{d_z} z(x)=z(x)$, we
see that $\alpha^{d_z}=1$. This yields that $d_z$ is the order of
element $\alpha $ in $GF(q)^{*}$. Consequently, $d_z|q-1$, $1 \le
d_z \le q-1$.

Finally, under the condition $gcd(n,q)=1$ the polynomial $x^{n}-1$
has~ $n$ different roots in $GF(q^h)$ field, where $h$~is the
multiplicative order of~$q$ modulo $n$~\cite{LID}. Denote by
$E(n)$ the multiplicative group of $n$-~th roots of unity over
$GF(q)$. Let $\xi \in E(n)$ be some~ $n$-~th~primitive root of
unity. Then the following set of elements $ \xi^0 = 1, \xi , \xi
^2, \ldots, \xi ^{n_z -1}, \ldots , \xi^{n-1}$ represents the
group $E(n)$. Since $n_z|n$, we have $E(n_z)\subset E(n)$, where
$E(n_z)$ is the multiplicative group of $n_z$-th roots of unity.
Moreover, taking into account the isomorphism of the groups,
having the same order~ \cite{KUR}, we can state that the subgroup
$\{ \alpha\}$, $\{ \alpha \} \subset GF(q)^{*}$, of order~$d_z$
belongs to ~$E(n_z)$ because $d_z|n_z$.

As mentioned above, $n_z$ is the period of $z(x)$, so that $n_z$
is the smallest divisor of~ $n$ such that the following congruence
$x^{r_z}\equiv \alpha \, {mod}(x^{n_{z}}-1)$ takes place.
Hence~$E(n_z)$ is the smallest group of ~$n$-~th roots of unity,
which contains $\{\alpha \}$. Since $\alpha = \xi^{r_z}$, we see
that the following elements $\xi^{r_{z}d_{z}}= 1, \xi^{r_z},
\xi^{2r_z}, \ldots,\xi^{(d_{z}-1)r_z}$ represent the subgroup~ $\{
\alpha \}$ in the group~$E(n_z)$. Besides that, the~decomposition
of~$E(n_z)$ relative to the subgroup~$\{ \alpha \}$ consists of
$r_z$ different cosets. Thus, the number $r_z$ is the index of
subgroup $\{ \alpha \}$ in the group of roots of unity, having the
smallest possible order. The theorem is proved.

{\bf Remark $2$.} Notice that when a parity-check polynomial of
code is some primitive polynomial of degree~$m$ and of
order~$n=q^{m}-1$, $m>1$, then  $r =R={(q^{m}-1)/(q-1)}$ and
$d=q-1$.~(see~\cite{MAKS}, \cite{PET}).

{\bf Corrollary $1$.} {\it The period} $n_z$, $n_z=r_z\cdot
d_z$,{\it of element} $z(x)$, ${z(x)\in J}$, {\it equals} $d_z$,
$d_z|(q-1)$, $1\le d_z \le (q-1)$,{\it if and only if the cycle}
$\{z(x)\}$ {\it is contained in one class of proportionality.}

{\bf Corollary $2$}. {\it All code words of any cyclic code,
havinq some length} $n$ {\it over} $ GF(q)$, {\it fall into some
equal-weight subsets and every such subset includes all
proportional to each other cycles.}

Besides that, consider some minimal ideal $J$, $J \subset  A_n$,
having an irreducible nonprimitive parity-check polynomial $h(x)$
of degree $m$ and of order $n$, $n \ne q^{m} - 1$, i.e., some
irreducible code of nonprimitive length.

{\bf Remark $3$.} The degree $m$ of the polynomial $h(x)$
coincides with the multiplicative order $h$ of the number $q$
modulo~$n$~\cite{LID}. Also, the order $n$ of the polynomial
$h(x)$ is some divisor of $q^{h}-1$. This means that the order $n$
can change in the following limits $1<n<q^{h}-1$, ${n \ne q-1}$.
If ${1< n < q-1}$, i.e., ${n\ne q-1}$, then some minimal
ideal~$J$, ${J\subset A_n}$, of dimention one contains only one
nonzero class of proportyonality. Consequently, $n=d$, $1< d <
q-1$, $d|q-1$.

{\bf Theorem $2$}. {\it Any cycle of minimal ideal} $J$, $
J\subset A_n$, {\it having some parity-check polynomial} $h(x)$
{\it of degree} $m$, $m>1$, {\it and of order} $n$, $n \ne
q^{m}-1$, {\it is contained in} $r$ {\it proportionality classes},
$1\le r \le R$, $r|R$, $R=(q^{m}-1)/(q-1)$. {\it Every such class
consists of} $d$, $1\le d \le (q-1)$, {\it different vectors of
cycle}, {\it and}

$$
                       d=gcd(q-1,n), \quad r=n/d.
\eqno(9)
$$
{\bf Proof}. All elements of minimal ideal~$J$, $J\subset A_n$,
have the same order $n$, $n=ord(h(x))$. Applying the theorem $1$
to some element $f(x)$, $f(x)\in J$, we have $n=r_{f}\cdot d_f$,
$d_f|(q-1)$, $1\le r_f\le n$, $1\le d_f \le (q-1)$, and also the
following equality
$$
               x^{r_f}\cdot f(x)= \gamma \cdot f(x),        \eqno (10)
$$
where $\gamma$ is some element of $GF(q)^*$. Evidently, if either
$r_f=1$, i.e., $d_f=n$, or $r_f=n$ and $d_f=1$, then equalities
(9) take place.  Hence, below we suppose that $1<r_f<n$,
$1<d_f<q-1$, and therefore $q-1<n<q^{m}-1$.

Taking into account $(9)$, we have $d_f \le d$. Let us show that
the strong unequality $d_f<d$ is impossible. Indeed, if $d_f <d$,
then the subgroup $\{\gamma \}$, where $\gamma$ is the element
from equality~$(10)$, belongs to some group of $n$-th roots of
unity, having the order $d$, because $d_f|d$. Since $d|(q-1)$, we
see that the subgroup $\{ \gamma \}$ belongs to $GF(q)^*$. This
means that the cycle $\{f(x)\}$ is contained in one class of
proportionality, i. e., $r_f=1$. But this fact contradicts to the
condition $r_f>1$. This implies that the strong unequality $d_f<d$
is impossible. Hence $d_f=d$. Because of an arbitrary choice of
$f(x)$ we can conclude that equalities~$(9)$ take place for any
element of~$J$. The theorem is proved~\cite{MLL2}.

{\bf Remark $4$.} It is necessary to note that the theorem $2$ is
valid only for irreducible codes of non--primitive length except
Reed-Solomon codes of length $n=q-1$ as it was shown above. In the
case of irreducible codes of primitive length  $n$, $n=q^{m}-1$,
$m>1$, the theorem~$2$ will be valid if and only if
$gcd(m,q-1)=1$. Indeed, when the last condition takes place, then
${gcd((q^{m}-1)/(q-1),q-1)=1}$~\cite{PET}. Thus,
$d=q-1={gcd(q^{m}-1,q-1)}$ that is $d=gcd(n,q-1)$. It follows that
the theorem~$2$ holds.

{\bf Remark $5$.} Notice that under the condition~${gcd(m,q-1)=1}$
the number $r$ from $(9)$ has no divisors of $(q-1)$ except~$1$,
so~${gcd(r,q-1)=1}$. This means that $gcd(r,d)=1$. Hence,
considering the fact that $r|R$ and also, taking into account
$(9)$ and the following equality $s=R(q-1)/n$ , we have
$r=gcd(R,n)$. Besides that, if either one from the two numbers $n$
and $q-1$ does not contain multiple prime divisors or the same
prime divisors of these numbers have the same degrees under the
decomposition of both $n$ and $q-1$, then the following equalities
$r=gcd(R,n)$, $gcd(r,d)=1$ also take place.

{\bf Corollary $3$.} {\it Both the number} $r$ {\it and the
number} $d$ {\it are the same numbers of all irreducible divisors
of polynomial} $x^{n}-1$ {\it over} $GF(q)$, {\it having the same
order}.

{\bf Corollary $4$.} {\it The number $R$, $R=(q^{m}-1)/(q-1)$, of
proportionality classes},  {\it of some irreducible code} $K$,{\it
having some length} $n$, $n\ne(q^{m}-1)$, $n=d\cdot r$, {\it over}
$GF(q)$ {\it field, consists of some} $v$ {\it different subsets.
And every such subset contains} $r$ {\it equal-weight
proportionality classes, i.e.,} $R={v\cdot r}$. {\it Besides that,
every subset includes} $b$ {\it equal-weight cycles}, $b=(q-1)/d$,
$1<b\le(q-1)$. {\it So that the number of all cycles for} $K$ {\it
equals} $s=v\cdot b$ {\it and} $gcd(r,b)=1$.

{\bf Proof.} According to the theorem $2$ any cycle of code $K$ is
contained in $r$, $r|R$, proportionality classes. Therefore the
number $v=R/r$ gives us the common quantity of different subsets
of $J$, each of which consists of $r$ classes, i.e., $R=v\cdot r$.
The number $b$, $b=(q-1)/d$, is the number of all different
equal-weight cycles,  contained in every such subset,
 which consists of some $r$ classes. Hence the number of all cycles
for $K$ is equal to $s=v\cdot b$. Since $d=gcd(q-1,n)$ we see that
$gcd(r,b)=1$. Actually, assuming the inverse, we would have been
able to decrease the number $d$, but it's impossible. The
corollary is proved.

{\bf Corollary $5$.} {\it The irreducible nonprimitive code} $K$
{\it is some equidistant code if} $s=b$.{\it Besides that, the
last equation is equivalent to the following ones}: $r=R$ {\it or}
$gcd(s,R)=1$.

{\bf Remark $6$.} Note that the condition $s=b$ was obtained in~
\cite{CLR,OYT}, but only for some subclass of irreducible
nonprimitive codes and under the following additional restriction
$gcd(b,m)=1$.

{\bf Corollary $6$.} {\it The weight of any element, belonging to
some irreducible nonprimitive code} $K$ {\it of length} $n$ {\it
over} $GF(q)$, {\it is multiple of the number} $d$,
$d=gcd(q-1,n)$.

{\bf Proof.} The weight of any element $z(x)$, $z(x) \in K$, of
order $n$, $n=ord(h(x))$, is equal to the number of such $j$, $0
\le j \le n-1$, for which the polynomial ${x^{j}\cdot z(x)}$ has
the following degree $n-1$. According to the theorem~$2$, the
number of such polynomials for the cycle ${z(x)}$, having degree
$n-1$, is equal to $w_{r}\cdot d$, where $w_{r}$ is the number of
polynomials, having degree $n-1$, among the first $r$ cyclic
shifts of $z(x)$, and $d= gcd(q-1,n)$. The corollary is proved.

In addition, consider some ideal $J$, $J\subset A_n$, of the form
$(2)$, having the parity-check polynomial $h(x)$ in terms of
$(3)$.

{\bf Theorem $3$}.{\it If the following  condition}
$gcd(h,q-1)=1$, {\it where} $h$~{\it is the multiplicative order
of number} $q~\mod~n$, {\it takes place}, {\it then any cycle of
set} $C$, $C\subset J$, {\it having some minimal polynomial of the
form} $(4)$, {\it is contained in} $r_c$ {\it proportionality
classes,} $r_c|R_c$, {\it and every such class includes} $d_c$,
$d_c|q-1$, {\it elements of cycle, that is} $n_c=r_{c}\cdot d_c$,
$n_c={ord}(c(x))$, {\it where} $R_c$~ {\it is the number of all
proportionality classes of set} $C$, {\it and}
$$
              r_c=gcd(R_c,n_c), \quad d_c=gcd(q-1,n_c).      \eqno (11)
$$
{\bf Proof}. It is sufficient to consider the case $k=2$ because
the general case can be obtained by the induction. Thus assume
that $c(x)=h_{1}(x)\cdot h_{2}(x)$, where $h_{i}(x)$ is of degree
$m_i$ and of order $n_i$, $n_i=q^{m_i}-1$, $1\le i \le 2$, is a
certain prime miltiplier of~$c(x)$. It is {known~\cite{LID}}, that
the number $m_i$, $1\le i \le 2$, equals either $h$ or some
divisor of this number. Hence, taking into account the
theorem~$2$, and also remarks $4$ and $5$, we have $n_i=r_i\cdot
d_i$, where $r_i=gcd(R_i,n_i)$, $d_i=gcd(q-1,n_i)$,
$gcd(r_i,q-1)=1$, $1\le r_i\le R_i$, $1\le d_i\le q-1$, and
$R_i={(q^{m_i}-1)/(q-1)}$, where $R_i$~is the number of all
proportionality classes of minimal ideal~$J_i$, $1\le i\le 2$.
Therefore the order~$n_c$, $n_c =lcm(n_1,n_2)= n_1\cdot
n_2/gcd(n_1,n_2)$ of the polynomial $c(x)$ can be rewritten as
$$
    n_c=r_{1}~d_1\cdot r_{2}~d_2/gcd(r_{1}d_1\cdot r_{2}d_2).   \eqno
(12)
$$
Since $gcd(r_i,q-1)=1$, we have $gcd(r_i,d_{1}\cdot d_2)=1$,
${1\le i \le 2}$. Thus $gcd(r_{1}r_2,d_{1}d_2)=1$. Hence
$gcd(gcd(r_{1},r_2),gcd(d_{1},d_2))=1$ so that $(12)$ may be
represented in the following form
$$
    n_c=(r_1~r_2/gcd(r_1,r_2))\cdot (d_1~d_2/gcd(d_1,d_2)).     \eqno
(13)
$$
Thus, $n_c=lcm(r_1,r_2)\cdot lcm(d_1,d_2)$. Now by $r_c$ and $d_c$
we denote $lcm(r_{1},r_2)$ and $lcm(d_{1},d_2)$, respectively.
Thus $n_c=r_c\cdot d_c$ and $gcd(r_c,d_c)=1$. Since $d_c|q-1$ and
$gcd(r_c,q-1)=1$, we obtain $d_c=gcd(q-1,n_c)$. Considering $(5)$,
it follows that $n_c|(R_c\cdot(q-1))$. Hence we have
$r_c=gcd(R_c,n_c)$. Consequently  the order $n_c$ of any element
of set $C$ is equal to the product of two relatively prime
numbers, i.e., $n_c=r_c\cdot d_c$, where $r_c=lcm(r_1,r_2)=gcd
(R_c,n_c)$ and $d_c=lcm(d_1,d_2)=gcd(q-1,n_c)$.

Furthermore, applying the theorem~$1$ to some element $a(x)\in C$
of period $n_c=r_a\cdot d_a$, we have
$$
       x^{r_a}\cdot a(x)=\theta a(x),                         \eqno
(14)
$$
where $\theta \in GF(q)^*$ is some element of order $d_a$, and
$r_a$ is the smallest natural number such that equality~$(14)$
takes place. Notice that the subgroup $\{ \theta \}$ has the order
$d_a$, ${d_a<d_c}$. If $d_c=1$, then $d_a=1$ and
$n_c=r_a=r_c=gcd(R_c,n_c)$, so that equalities~$(11)$ hold. For
this reason below we suppose that $d_c>1$. If under this condition
the number $r_c$ is equal to one, then $n_c=d_c$=$gcd(q-1,n_c)$
and the theorem is valid. Therefore below we suppose that both
$r_c>1$ and $d_c>1$.

Evidently, $d_a\le d_c$. Now let us show that the
inequality~$d_a<d_c$ is not possible. Indeed, if $d_a<d_c$, then
we come to the following conclusion. The subgroup $\{ \theta \}$,
where $\theta$~is the element from equality~$(14)$, belongs to
some subgroup of $GF(q)^*$, having the order~$d_c$, because
$d_a|d_c$. Since $d_c$ is some divisor of $n_c$, then, considering
the uniqueness of subgroups, having the same order, the subgroup
$\{ \theta \}$ belongs to some group of $d_{c}$-th roots of unity.
This implies that the smallest group of roots of unit, containing
$\{ \theta \}$, has an order, which either less or equals $d_c$.
Thus, both the period of $a(x)$ and the order of ~$c(x)$ must be
either less or equal to $d_c$. This yields that the order
of~$c(x)$ must be some divisor of $q-1$. But this fact contradicts
the condition $r_c>1$.  Hence the inequality $d_a<d_c$ is
impossible so that $d_a=d_c$ and $r_a=r_c$. Because of an
arbitrary choice of $a(x)$ equalities~$(11)$ take place for any
element of set $C$. The theorem is proved.

{\bf Corollary $7$.} {\it The order} $n_c$ {\it of reducible
factor} $c(x)$ {\it of polynomial} $x^{n}-1$, {\it having some
degree}~$m$ {\it over} $GF(q)$, {\it is some divisor of number}
$q^{m}-1$, {\it if} ${gcd(h,q-1)=1}$, {\it where} $h$ {\it
is~the~multiplicative~order~of number}~$q~\mod~n$.

{\bf Remark $7$.} In terms of condition ${gcd(h,q-1)=1}$, where
$h$ is the multiplicative order of number $q~\mod~n$, the theorem
$3$ is valid for cyclic codes of both the primitive and the
nonprimitive length. Also, taking into consideration the remark
$5$, the order of any reducible factor of the polynomial
$(x^{n}-1)$ over $GF(q)$ of degree $m$, is some divisor of the
number $(q^{m}-1)$, if {\it gcd}${(h,q-1)=1}$.

{\bf Theorem $4$.} {\it Any cycle of set} $C$, $C\subset J$, {\it
having some minimal polynomial} $c(x)$ {\it of the type} $(4)$
{\it and of order} $n_c =lcm(n_1,n_2,\ldots, n_k)$, {\it where}
$n_i\ne q^{m_i}-1$, $1\le i\le k$, {\it is contained in} $r_c$,
$r_c|R_c$, $1\le r_c\le R_c$, {\it classes of proportionality and
every such class includes} $d_c$, $1\le d_c\le q-1$, {\it elements
of cycle, where} $R_c$ {\it is the number of all proportionality
classes of} $C$, {\it and}
$$
                 d_c=gcd(q-1,n_c), \quad r_c=n_c/d_c.
\eqno (15)
$$
{\bf Proof.} It is sufficient to assume that $k=2$ because the
general case can be obtained by the induction. This implies that
$c(x)=h_{1}(x)h_{2}(x)$, where $h_{i}(x)$ is some prime divisor of
equality $(3)$, having some degree~$m_i$, and of order~$n_i$,
${n_i\ne q^{m_i}-1}$, $1\le i \le 2$. This yields that the
theorem~$2$ holds for the minimal ideal $J_i$, $1\le i \le 2$.

Evidently, if one of the two numbers, i.e., either $d_c$ or $r_c$
is equal to one, then equalities~$(15)$ hold. For this reason
below we assume that both $r_c>1$ and $d_c>1$.

Let $z(x)$ be some vector of set $C$. According to the theorem~$1$
some cycle $\{ z(x) \} \subset C$ is contained in $r_z$ classes
and every such class includes $d_z$ different elements of this
cycle. Thus $n_c=r_z~d_z$, $1\le d_z \le q-1$, $d_z|q-1$, $1\le
r_z \le n_c$. And in addition, the following equality takes place
$$
              x^{r_z}\cdot z(x) = \beta~z(x),                   \eqno
(16)
$$
where $\beta$ is some element of $GF(q)^*$ and the subgroup $\{
\beta \}$, $\{ \beta \} \subset GF(q)^*$, has the order $d_z$.
Evidently, $d_z \le d_c$. Let us show that the number $d_z$ can
not be smaller than $d_c$. Assume the inverse, i.e. let $d_z$ be
less than $d_c$. Since at least one of the two numbers ether $r_1$
or $r_2$ is not equal to one, we see that at least one of the
numbers $n_i$, $1\le i\le 2$, is more than $q -1$, as was
established in the theorem 2. This means that
$$
               n_c=lcm(n_1,n_2)> (q - 1).                       \eqno
(17)
$$
Further, since $d_z|d_c$, we see that the subgroup~ $\{ \beta \} $
of order $d_z$ belongs to some subgroup of $GF(q)^*$, having the
order $d_c$, where $\beta$ is the element from equality~$(16)$.
Due to  the uniqueness of groups, having the same order, the
subgroup~ $\{\beta \}$ belongs to the group of $d_c$-th roots of
unity  because $d_c|n_c$. It follows that the smallest group of
the roots of unity, which contains the subgroup $\{ \beta \}$, has
the order less or equals to $d_c$. This implies that both the
period of $z(x)$ and the order of $c(x)$ is some divisor of $q-1$.
But this fact contradicts to $(17)$. It follows that our
assumption is not true, so that $d_z=d_c$, and $r_z = r_c$.

Besides that, since the number $d_c$  is the same number for every
cycle of set~$C$, we see that every subset, consisting of $r_c$
proportionality classes, contains the same number of cycles, which
is equal to $b_c=q-1/d_c$. Moreover, since $d_c=gcd(q-1,n_c)$, we
obtain $gcd(b_c,r_c)=1$. Hence, taking into account the following
equality $s_c=R_c(q-1)/n_c$, which follows from $(5)$, we  see
that $r_c$ is some divisor of $R_c$ and $s_c=v_{c}b_c$, where
$v_c$ is equal to $R_c/r_c$. The theorem is proved.

{\bf Corollary $7$.} {\it (Equidistance signs of the subset} $C$,
$C\subset J$, {\it having some minimal polynomial of the form}
$(4)$.)

{\it All vectors of the subset} $C$, $C\subset J$, {\it having
some minimal polynomial} $c(x)$ {\it of order~} $n_c=lcm(n_1,
n_2,..., n_k)$, $n_i\ne q^{m_i}-1$, $1\le i\le k$, $k>1$, {\it
have the same weight if at least one of the following conditions
holds:} $1.$ $s_c=b_c$, $1< b_{c} \le q-1$, $2.$ $r_c=R_c$, $3.$
$gcd(s_c,R_c)=1$.

{\bf Corollary $8$.} {\it The order of the reducible factor}~
$c(x)$ {\it of the polynomial} $x^{n}-1$ {\it over} $GF(q)$, {\it
having some degree} $m$, {\it is some divisor of the number}
$q^{m}-1$, {\it if the decomposition of the polynomial} $c(x)$
{\it into prime multiples does not contain any primitive
polynomial.}

{\bf Remark $8$.} Note that corollaries $(6)$ and $(8)$ show us in
what cases corollary~$3.4$ from~\cite{LID} takes place for some
reducible polynomial of the degree $m$ over $GF(q)$.

Also, consider some cyclic code, having the following parity-check
polynomial
$$
            h(x)= \prod_{i=1}^{t} h_{i}(x),                  \eqno (18)
$$
where $h_i(x)$ is an irreducible polynomial over $GF(q)$, $q>2$,
of  degree $m_i$ and of  order $n_i=(q^{m_i}-1)/b_i$, $b_i|q-1$,
$1\le b_i \le q-1$, $1\le i \le t$, $gcd(m_i,m_j)=1$ and
$gcd(b_i,b_j)=1$, provided $i\ne j$, $1\le i,j \le t$, so that the
order $n$ of $h(x)$ equals $n=lcm(n_1,n_2, ... ,n_t)$.

It is worth mentioning that in \cite{CLR} and \cite{OYT} the
following cases of polynomial $(18)$ have been considered. Namely
$t=1$ and $t=2$, and besides that, with some additional
restrictions, which can be omitted. Also, some particular case of
polynomial $(18)$, that is provided $b_i=1$ for all $i$, $1\le i
\le t$, was obtained in \cite{TEN}. But there are  some
unnecessary restrictions in this paper too. Also, there are some
mistakes in that paper. Namely, the order of the product for some
two polynomials from~$(18)$ was defined incorrectly in~
\cite{TEN}.

Finally, using the results obtained above, we have found the
minimal distance of code, having the parity-check polynomial
$(18)$,~(see \cite {MLL3}).

Denote by $M$ the set of degrees for the polynomials $h_i(x)$,
$1\le i \le t$, from the eq. $(18)$ that is $M=(m_1,m_2, ...
,m_t)$. Let the number $m$ denote the degree of polynomial $h(x)$,
i. e., $m = m_1 + m_2 + ... + m_t$. By $M_{j,k}$, $1\le j \le
C_{t}^{k}$, $1 \le k \le t-1$, denote $j^{,}s$ $k$-subset of $M$,
where $C_{t}^{k}$ is the binomial coefficient.

Thus $M_{j,k}=(m_{j_{1}}, m_{j_{2}}, ... , m_{j_{k}})$, $1\le j
\le C_{t}^{k}$, $1\le k \le t-1$. At last, denote by $m_{j,k}$ the
sum of degrees, from the subset $M_{j,k}$, so that
$m_{j,k}=m_{j_{1}} + m_{j_{2}} + ... + m_{j_{k}}$, $1\le j \le
C_{t}^{k}$, $1\le k \le t-1$. It is obvious that $1\le m_{j_{k}}<
m$. Let us remark that provided $k=1$ the number $m_{j,k}=m_{j,1}=
m_{j}$ and $1\le j \le t$ because $C_{t}^{1}=t$.

{\bf Theorem $5$}.{\it The minimal distance of code, having the
parity-check polynomial}~$(18)$, {\it has the following form}
$$
  d_{min}=q^{m-1}-\sum_{j=1}^{C_{t}^{1}} q^{m_{j,1}-1}-
         - \sum_{j=1}^{C_{t}^{2}} q^{m_{j,2}-1} - ... -
          - \sum_{j=1}^{C_{t}^{t-1}} q^{m_{j,t-1}-1},     t\ge 2,
$$

$$
                 d_{min}=q^{m-1}(q-1)/b,            t=1, 1\le b \le
q-1,
$$

$$
                  d_{min}=q^{m-1}(q-1),                   t=1, b=1.
$$
In conclusion I should like to express my sincere gratitude to
{L.A.Bassaligo}, M.I. Boguslavskii and  E.T. Akhmedov  for helping
me to correct some mistakes in the original version of my paper.

\end{document}